\documentclass[runningheads]{llncs}

\usepackage[T1]{fontenc}
\usepackage{graphicx}
\usepackage{todonotes}
\usepackage{hyperref}
\usepackage{mathtools}
\usepackage{amsfonts}
\usepackage{subcaption}
\usepackage[ruled,longend]{algorithm2e}
\usepackage[capitalize]{cleveref}

\newcommand \dd{\,\mathrm{d}}

\newcommand \Hs{\mathcal{H}^{d-1}}
\newcommand \Om{\Omega}

\newcommand \eps{\varepsilon}
\newcommand \1{\mathbf{1}}

\newcommand \Per{\text{Per}}

\newcommand \N{\mathbb{N}}

\newcommand \R{\mathbb{R}}

\renewcommand \R{\mathbb{R}}

\newcommand{\id}{\operatorname{id}}

\DeclareMathOperator{\dist}{dist}
\renewcommand{\d}{\operatorname{d}}
\DeclareMathOperator{\sign}{sign}
\newcommand{\dive}{\operatorname{div}}

\newcommand{\Vol}{\operatorname{Vol}}
\DeclareMathOperator{\argmin}{arg\,min}

\begin{document}
\title{Meshless Shape Optimization using\\Neural Networks and Partial Differential Equations on Graphs}
\titlerunning{Meshless Shape Optimization}
%
\author{
Eloi Martinet\inst{1} 
\and
Leon Bungert\inst{1,2}
}
\authorrunning{E. Martinet, L. Bungert}
%
\institute{Institute of Mathematics, University of Würzburg, Germany
\and
Center for Artificial Intelligence and Data Science (CAIDAS),\\University of Würzburg, Germany\\
\email{\{eloi.martinet,leon.bungert\}@uni-wuerzburg.de}}
\maketitle              
\begin{abstract}
Shape optimization involves the minimization of a cost function defined over a set of shapes, often governed by a partial differential equation (PDE). 
In the absence of closed-form solutions, one relies on numerical methods to approximate the solution. The level set method---when coupled with the finite element method---is one of the most versatile numerical shape optimization approaches but still suffers from the limitations of most mesh-based methods.
In this work, we present a fully meshless level set framework that leverages neural networks to parameterize the level set function and employs the graph Laplacian to approximate the underlying PDE. Our approach enables precise computations of geometric quantities such as surface normals and curvature, and allows tackling optimization problems within the class of convex shapes.

\keywords{Shape Optimization \and Neural Networks \and Graph Laplacian.}
\end{abstract}

\section{Introduction}

Shape optimization addresses problems of the form 
\begin{equation}
    \label{eq:generic_so}
    \min_{\Om \in \mathcal A} J(\Om)
\end{equation}
where $\mathcal A$ is a set of shapes in $\R^d$ satisfying some geometrical constraints and $J: \mathcal{A} \to \R$ is a cost function. In the majority of shape optimization problems, the function $J$ depends on the solution of a partial differential equation (PDE) posed on the domain. In most cases including industrial applications, an exact solution to \labelcref{eq:generic_so} cannot be found and a numerical approximation must be sought.
Among the many different numerical shape optimization methods, one of the most versatile approaches is the level set method~\cite{Dapogny2023}. Contrary to classical methods based on mesh deformation, the level set method represents shapes $\Omega\in\mathcal{A}$ as the $0$-sublevel set of a function, which has the advantage to allow for topology changes in the course of the optimization process. 
However, it also introduces new difficulties, like the need of an \textit{ersatz material}, i.e., an artificial material outside of the shape $\Omega$. Depending on the PDE, deriving an appropriate ersatz material formulation can be challenging~\cite{martinet_numerical_2023}. 
Other issues are the steepening of the level set function which requires redistancing \cite{allaire2021shape} or the fact that the level of detail of the optimal shape is determined by the quality of the underlying mesh. Some of these issues can be alleviated by a re-meshed level set approach \cite{Allaire2014Dec} which, however, can be costly since it requires to re-mesh the domain and reinitialize the level set function at each iteration.

\subsection{Our contribution}

\begin{figure}[!t]
    \centering
    \begin{subcaptionblock}{0.32\textwidth}
        \includegraphics[width=\textwidth]{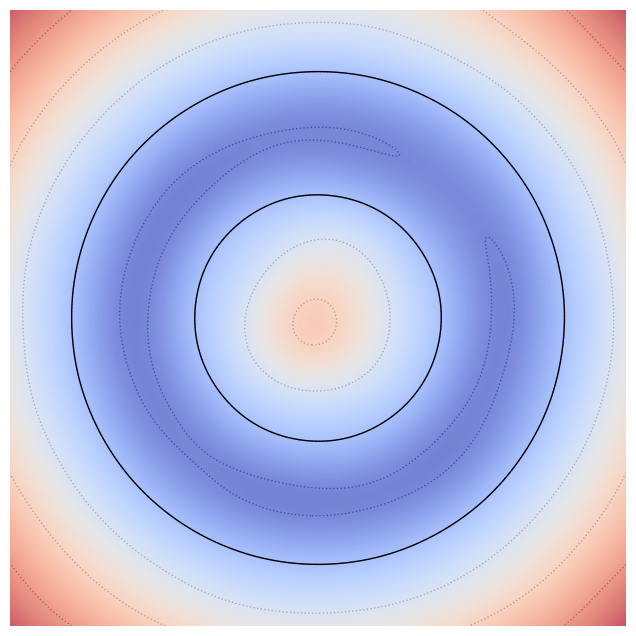}
        \caption{Level set function}
    \end{subcaptionblock}
    \begin{subcaptionblock}{0.32\textwidth}
        \includegraphics[width=\textwidth]{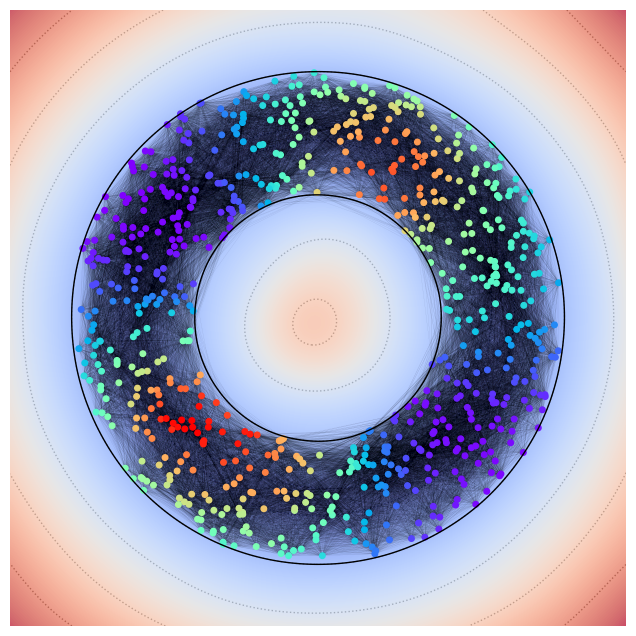}
        \caption{Solution of graph PDE}
    \end{subcaptionblock}
    \begin{subcaptionblock}{0.32\textwidth}
        \includegraphics[width=\textwidth]{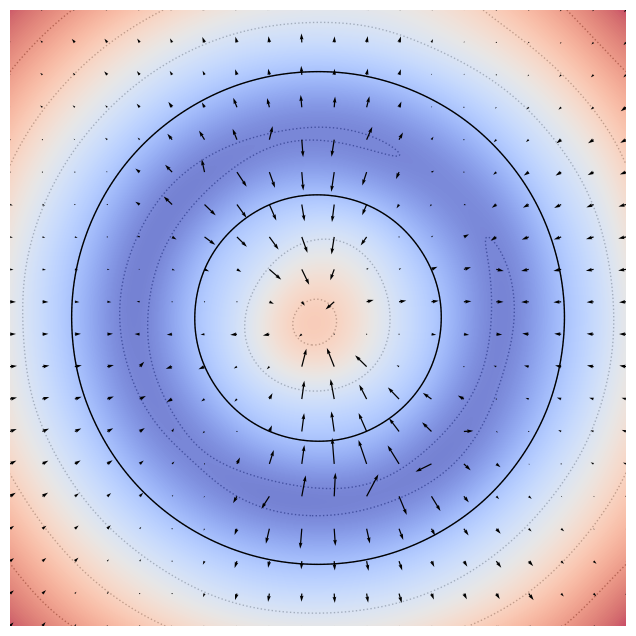}
        \caption{Shape derivative}
    \end{subcaptionblock}
    \caption{
        \label{fig:method}
        Key elements of the proposed method.
    }
\end{figure}

In this paper, we introduce a new level set method that is completely meshfree and does not require any ersatz material. Our method is based on a neural network representation of the level set function, while the PDE is solved using a graph Laplacian approach, see \cref{fig:method}.
We show that this approach allows us to naturally take into account some geometrical constraints as convexity, and seamlessly adapts to three-dimensional shape optimization problems, where meshing can be costly. Moreover, common geometric quantities like normals and curvature can be easily and precisely computed using automatic differentiation. 
Finally, by training the neural network such that it approximates the signed distance function of the shape, no redistancing is required.
We will demonstrate the potential of our method on three different shape optimization problem: the maximization of eigenvalues of the Neumann Laplacian under volume constraint, a Poisson problem with Dirichlet boundary conditions, and the minimization of eigenvalues of the Neumann Laplacian under perimeter and convexity constraints.

\subsection{Related work}

Several meshless shape optimization methods have already been studied in the literature. 
For instance, one can represent the boundary of a 2D shape by a truncated Fourier series, or a decomposition on the basis of (hyper)spherical harmonics in higher dimensions.
The PDE is then solved using some meshless method, such as the method of fundamental solutions \cite{Bogosel2016Nov,antunes_numerical_2017}. 
Although computationally efficient, this method can only deal with very specific equations on simply connected domains.
Previous works like \cite{deng_parametric_nodate} already used a neural network representation of the level set function, which falls in the broader category of \textit{parametrized level set methods} \cite{Cui2021Apr}. However, to the knowledge of the authors, the various parametrized level set methods for shape optimization still rely on the finite element method to compute the solution of the underlying PDE, and hence suffer from some of the previously detailed meshing issues and from the need of an ersatz material. Moreover, the flexibility of a neural network representation of the level set function for computing geometric quantities and taking geometric constraints into account has not been yet studied.

In recent years, physics-informed neural networks (PINNs) emerged as a new meshless method for PDEs \cite{raissi2017physics}. 
In this context the solution of a PDE is represented by a neural network, the parameters of which are optimized such that the PDE is solved approximately on a number of collocation points. 
The advantage of this approach is that all gradients involved in optimization can be automatically computed using backpropagation.
The main drawbacks are the lack of general convergence results of PINNs toward the solution of the PDE and the computational cost compared to well-established methods such as the finite element method \cite{grossmann2023physicsinformedneuralnetworksbeat}.
In the context of shape optimization, PINNs were used, for example, in \cite{Belieres--Frendo2024Jul}, where the shape is represented as a symplectic transformation of a given domain and the associated symplectomorphism as a neural network. 
While being meshless, this method does not allow for topology changes. In the context of linear elasticity, some fully meshless methods implemented the \emph{solid isotropic material with penalization} approach using only neural networks \cite{Zehnder2021Dec}. However, this method relaxes the original shape optimization problem into a density optimization problem, and does not readily generalize to other PDEs.

\section{Description of the method}

In this section, we detail every step of our meshless shape optimization algorithm. A self-contained version in algorithmic form can be found in \cref{alg:mso} at the end of this section.

\subsection{Shape derivative and level set representation}

Numerical optimization of \labelcref{eq:generic_so} usually relies on a gradient descent algorithm which computes a sequence of shapes $(\Om_k)_{k \in \N}$ starting at $\Omega_0\subset\R^d$ according to 
\begin{equation}
    \label{eq:gradient_descent}
    \Om_{k+1} = (\id + \tau_k V_k)(\Om_k).
\end{equation}
Here, $\tau_k> 0$ is a step size at the $k^\text{th}$ iteration and $V_k :\R^d\to\R^d$ is a vector field chosen in a way that $J(\Om_{k+1}) \leq J(\Om_k)$ where $J$ is the cost function in \labelcref{eq:generic_so}. 
In order to compute such a vector field, one needs the concept of a \textit{shape derivative}. 
For a given vector field $V:\R^d \to\R^d$, the shape derivative of $J$ at $\Om$ in the direction $V$ is defined (when it exists) as 
\begin{equation}
    \label{eq:shape_der}
    J'(\Om,V) := \lim_{t \to 0} \frac{J\left((\id +tV)(\Om)\right) - J(\Om)}{t}.
\end{equation}
If we find some $V$ with $J'(\Om,V)<0$, we can decrease the value of $J$ by taking a small step in the direction $V$. Under mild assumptions, the Hadamard structure theorem \cite{henrot_variation_2005} asserts the existence of $f : \partial \Om \to \R$ such that
\begin{equation}
    \label{eq:hadamard}
    J'(\Om,V) = \int_{\partial \Om} f(V\cdot n)\dd\Hs,
\end{equation}
where $n$ is the normal on $\partial \Om$. Hence, taking $V := -f_\Om n_\Om$ (where $f_\Om$ and $n_\Om$ are respective extensions of $f$ and $n$ to $\R^d$) is a descent direction since $J'(\Om,V)<0$.

Let $D =[0,1]^d$ be the computational domain and $\Om \subset D$. In the level set method, the shape $\Om \subset D$ is represented by a function $\phi:D \to \R$ which verifies 
\begin{equation*}
    \phi(x) < 0\text{ for } x \in \Om,
    \qquad
    \phi(x) = 0\text{ for } x \in \partial \Om,
    \quad 
    \text{and}
    \quad 
    \phi(x) > 0\text{ for } x \in D \setminus \bar{\Om}.
\end{equation*}
In our approach, the level set function of a domain $\Om$ is represented by a neural network $\phi_\theta:D\to\R$ where $\theta$ denotes the parameters of the network. Computing~\labelcref{eq:gradient_descent} where $\Om_k = \{\phi_{\theta_k} < 0 \}$ then amounts to computing parameters $\theta_k$.

\subsection{Computation of geometric quantities}\label{subseq:geom}

One of the main advantages of representing of $\Om$ by a neural network $\phi_\theta$ is the easy and accurate computation of geometric quantities using automatic differentiation. 
For instance, the extended normal vector field $n_\Om$ of $\partial \Om$ is given~by
$$
    n_\Om(x) := \frac{\nabla \phi_\theta(x)}{|\nabla \phi_\theta(x)|}
$$
for all $x \in D$ with $\nabla\phi_\theta(x)\neq 0$, which can be automatically computed using backpropagation. Another quantity of interest in shape optimization is the \textit{mean curvature} of a set $\Om$, which appears for instance in the shape derivative of the perimeter functional and can be computed as
$$
    \kappa_\Om(x) = \dive n_\Om(x)  = \dive \left(\frac{\nabla \phi_\theta(x)}{|\nabla \phi_\theta(x)|}\right).
$$
Furthermore, using a Monte Carlo approach we approximate integral quantities like the volume $\Vol(\Om)$ and the perimeter 
$$
    \Per(\Om) := \int_{\partial \Om} n_\Om \cdot n_\Om \dd\Hs = \int_{\Om} \dive n_\Om \dd x = \int_{\Om} \kappa_\Om \dd x
$$
based on i.i.d. random variables $\{x_i\}_{1\leq i \leq N}$, distributed uniformly in $D$, as
$$
    \Vol(\Om) \approx \frac{1}{N} \sum_{i=1}^N  \1_{\{\phi_\theta < 0\}}(x_i) \quad \mbox{and} \quad \Per(\Om)  \approx \frac{1}{N} \sum_{i=1}^N  \1_{\{\phi_\theta < 0\}}(x_i) \kappa_\Om(x_i).
$$

\subsection{Advection of the level set}

While the advection of the level set usually amounts at solving the so-called level set advection equation \cite{Dapogny2023}, we chose to adopt a different approach which proved to be more accurate in our context while requiring equivalent computation time. For this we let $\d(x,A):=\inf_{y\in A}|x-y|$ denote the distance of a point $x\in\R^d$ to a subset $A\subset\R^d$ and let $\dist_\Om$ denote the signed distance function of $\Om$: 
\begin{equation*}
    \dist_\Om(x) := 
    \begin{dcases}
        -\d(x, \partial \Om) & \mbox{ for } x \in \overline\Om, \\ 
        \phantom{-}\d(x, \partial \Om) & \mbox{ for } x \in D \setminus \overline \Om. 
    \end{dcases}
\end{equation*}
For each time step $k$ in \labelcref{eq:gradient_descent} we aim at finding parameters $\theta_{k+1}$ such that 
$$
    \phi_{\theta_{k+1}}(x) \approx \dist_{\Om_{k+1}}(x)\qquad 
    \forall x\in D.
$$
Motivated by \labelcref{eq:gradient_descent}, we first sample the boundary $\partial \Om$ with points $\{y_j\}_{1 \leq j \leq N_y}$, and sample the domain $D$ with points  $\{x_i\}_{1 \leq i \leq N_x}$ which we then advect by
$$
    \tilde y_j := y_j + \tau_k V_k(y_j) \quad \mbox{and} \quad \tilde x_i := x_i + \tau_k V_k(x_i)
    \qquad\forall i,j.
$$
The points $\{\tilde y_j\}_{1 \leq j \leq N_y}$ are then a discretization of the next interface $\partial \Om_{k+1}$. We use them to compute the following approximation of $\dist_{\Om_{k+1}}$ at the points $\tilde x_i$:
\begin{align}\label{eq:proxy_distance}
    \widehat{\dist}_{\Om_{k+1}}(\tilde x_i) := \sign(\phi_{\theta_k}(x_i))\d(\tilde x_i, \{\tilde y_j\}_j). 
\end{align}
Indeed, the quantity $\d(., \{\tilde y_j\}_j)$ is an approximation of the (unsigned) distance function to the boundary $\partial \Om_{k+1}$, while $\sign(\phi_{\theta_k})$ allows us to recover the sign.
Finally, the network parameters $\theta_{k+1}$ are then computed as the minimizer of 
\begin{equation}\label{eq:signed_dist_fitting}
    \min_{\theta} \left\{ \frac{1}{N_x} \sum_{i=1}^{N_x} \left|\phi_\theta(\tilde x_i) - \widehat{\dist}_{\Om_{k+1}}(\tilde x_i) \right|^2 + \frac{\alpha}{N_y} \sum_{j=1}^{N_y} \left|\phi_\theta(\tilde y_j) \right|^2 \right\},
\end{equation}
where the parameter $\alpha\geq 0$ allows enforcing a greater accuracy on boundary points.

\begin{remark}
We note that there is an alternative way to approximate the signed distance function of $\Omega_{k+1}$ without the need to compute the distance $\widehat{\dist}_{\Omega_{k+1}}$ defined in \labelcref{eq:proxy_distance} for all points $\tilde x_i$.
Using the fact that
\begin{align*}
    \dist_A = 
    \argmin\left\lbrace
    -\int_{\R^d} \phi(x) \left(\frac12 - \1_A(x)\right)\d x\,:\,\operatorname{Lip}(\phi)\leq 1,\,\phi\vert_{\partial A}=0\right\rbrace,
\end{align*}
we can compute network parameters by letting $\phi_\theta$ be a $1$-Lipschitz neural network~\cite{anil2019sorting,hasannasab2020parseval}, replacing the first term in \labelcref{eq:signed_dist_fitting} by $-\frac{1}{N_x}\sum_{i=1}^{N_x}    \phi_\theta(\tilde x_i)\sign(\phi_{\theta_k}(x_i))$, and choosing a large value of $\alpha>0$.
However, in our numerical experiments this approach---although undeniably elegant---turned out to be less accurate than \labelcref{eq:signed_dist_fitting} which is why we only report results using \labelcref{eq:signed_dist_fitting}.
\end{remark}

\subsection{Sampling of the interface}

The previous advection procedure makes it necessary to uniformly sample the boundary $\partial \Om$ of $\Om \subset \R^d$ by a set of points $\{y_j\}_{1 \leq j \leq N}$. 
For this we solve
\begin{equation}
    \label{eq:riesz}
    \min_{y_1, \dots ,y_N \in \partial \Om} \sum_{1 \leq i\neq j \leq N} \frac{1}{\left|y_i-y_j\right|^s}
\end{equation}
with $s>d-1$. The optimal distribution of points $\{y_j^N\}_{1\leq j \leq N}$ can be seen as the equilibrium state of $N$ charged particle constrained on the surface $\partial \Om$. 
Note that the \textit{poppy-seed bagel theorem} \cite{borodachov_discrete_2019} states that under mild assumptions, the empirical measures $\nu_N := \frac{1}{N}\sum_{j=1}^N \delta_{y_j^N}$ converge weakly-$\star$ to the uniform Hausdorff measure on $\partial \Om$ as $N\to\infty
$.
In our case, the constraints $y_i \in \partial \Om$ are equivalent to $\phi_\theta(y_i)=0$ for all $i \in \{1,\dots,N\}$ and the resulting constrained optimization problem \labelcref{eq:riesz} can be numerically solved using IPOPT~\cite{ipopt}.


\subsection{Solving PDEs using the graph Laplacian}

As mentioned in the introduction, two of the shape optimization problems considered thereafter involve eigenvalues of the Neumann Laplacian. 
As was shown in \cite{calder_improved_2020} and the references therein, these can be approximated by the ones of the graph Laplacian on random geometric graphs. While recently the graph Laplacian has mainly been studied in the context of unsupervised \cite{von_luxburg_tutorial_2007} and semi-supervised \cite{zhu2003semi,calder2023rates,calder2020poisson,bungert2024convergencerates} learning, we will use it as a numerical PDE solver.

Let $\Om\subset \R^d$ as before and $\mu$ be a probability distribution supported on $\Om$ with a positive continuous density $\rho $. 
Suppose that $V^{(n)} = \{x_1,\dots,x_n\}$ are i.i.d. samples of $\mu$ and let $\eta:(0,\infty)\to [0,\infty)$ be compactly supported and non-increasing. 
We define the \textit{weight matrix} $W^{(n)}$ as
\begin{equation*}
    W_{ij}^{(n)} = \frac{2}{n\sigma_\eta\eps_n^{d+2}} \eta\left( \frac{|x_i - x_j|}{\eps_n} \right),
    \qquad 
    i,j=1,\dots,n,\,i\neq j,
\end{equation*}
where $\eps_n>0$ is a scaling parameter and $\sigma_\eta>0$ is a constant depending only on $\eta$ and $d$. 
The \textit{degree matrix} $D^{(n)}$ is defined as the diagonal matrix with entries $D^{(n)}_{ii} = \sum_{i=j}^{(n)} W^{(n)}_{ij}$.
The (unnormalized) \textit{graph Laplacian} $L^{(n)}$ is then defined as
$$
    L^{(n)} := D^{(n)} - W^{(n)}
$$
As a positive semi-definite symmetric matrix, the graph Laplacian has eigenvalues
$
    0 = \lambda_0^{(n)} \leq \lambda_1^{(n)} \leq \dots \leq \lambda_n^{(n)}.
$
With a suitable choice of $\eps_n$ going to zero at a certain rate in $n$, for every $k \in \mathbb{N}$ it holds almost surely \cite{calder_improved_2020} that
\begin{equation*}
    \lim_{n \to \infty}  \lambda_k^{(n)} = \tilde\lambda_k,
\end{equation*}
where $\tilde \lambda_k$ is the $k^\text{th}$ non-trivial eigenvalue of 
\begin{equation}
    \label{eq:neumann_density}
    \begin{dcases}
        -\frac{1}{\rho} \dive(\rho^2 \nabla u) = \tilde \lambda_k u &\mbox{ in } \Om,\\
        \partial_n u = 0 &\mbox{ on } \partial \Om.
    \end{dcases}
\end{equation}
In \cite{calder2022lipschitz} it was shown that also the corresponding eigenvectors of $L^{(n)}$ converge uniformly and in a Lipschitz-type semi-norm to the eigenfunctions of \labelcref{eq:neumann_density}, restricted to the vertices $V^{(n)}$.
Hence, letting $\bar u \in \R^n$ be an eigenvector associated to $\lambda_k^{(n)}$, we can approximate an eigenfunction $u$ associated to $\tilde\lambda_k$ by $u(x_i)\approx \bar u_i$.
Using Taylor's formula we can also derive an approximation of the gradient of~$u$:
\begin{equation}
    \label{eq:approx_gradient}
    \nabla u(x_i) \approx \frac{d}{D_{ii}^{(n)}} \sum_{j=1}^n W_{ij}^{(n)} \frac{\bar u_j - \bar u_i}{|x_j-x_i|^2}(x_j -x_i).
\end{equation}
The solution of the PDE and its gradient are usually used to derive the function~$f$ in \labelcref{eq:hadamard}. In this graph approach, $f$ will be a function defined over the vertices $V^{(n)}$. However, in order to define the advection vector field, we will need to define it over $\R^d$ by taking (for instance) its nearest neighbor interpolation $f_\text{NN}$.

Here is the complete algorithm:
\begin{algorithm}
\caption{\label{alg:mso} Meshless Shape Optimization Algorithm}
\KwIn{$k_{\rm max}, k_{\rm sample}, n, N_x, N_y,q \tau$}
\KwOut{Optimized shape and neural network parameters}

\tcc{Initialization}
Initialize the neural network $\phi_{\theta}$\;
Sample the boundary of $\{\phi_{\theta} < 0\}$ by points $y_1, \dots, y_{N_y}$ according to \labelcref{eq:riesz}\;

\For{$k \leftarrow 0$ \KwTo $k_{\rm max}$}{
    \tcc{Computation of the shape derivative}
    Set $\Omega\leftarrow \{\phi_{\theta} < 0\}$\;
    Compute $\Vol(\Omega_k)$ using Monte Carlo integration\;
    Draw $V^{(n)} = \{z_1, \dots, z_n\}$ as i.i.d. samples from $\Omega_k$\;
    Compute the graph Laplacian $L^{(n)}$ \;
    Solve the graph PDE \tcp{e.g, approximating \cref{eq:neumann_ev}}
    Compute the shape derivative $f$ on $V^{(n)}$ \tcp{see e.g. \cref{eq:sd_neumann}}
    \BlankLine
    \BlankLine
    
    \tcc{Advection of the level set}
    Sample $\{x_1, \dots, x_{N_x}\}$ uniformly from domain $D$\;
    \For{$i \leftarrow 1$ \KwTo $N_x$}{
        Compute normal vector $n_{\Omega}(x_i)$ using automatic differentiation\;
        Set velocity field $ V(x_i) \leftarrow f_\text{NN}(x_i) n_{\Omega}(x_i)$\;
        Update position: $\tilde x_i \leftarrow x_i + \tau V(x_i)$
    }
    \tcp{Same for $\{y_j\}_j$}

    Compute the next signed distance function \cref{eq:proxy_distance} \;
    Update neural network parameters according to \cref{eq:signed_dist_fitting}\;

    \BlankLine

    $\{y_j\}_j \leftarrow \{\tilde y_j\}_j$ \tcp{Update the interface points}
    \If{$k \bmod k_{\rm sample} = 0$}{
        Resample the boundary points\;
    }
}
\end{algorithm}

\section{Numerical examples}

\subsection{Optimization of Neumann Eigenvalues using the graph Laplacian}

Let $\Omega \subset \R^d$ be an open set with Lipschitz boundary. For $k \geq 1$, we consider the maximization problem
\begin{equation}
    \label{pb:neumann}
    \max_{\Vol(\Om)=1 } \mu_k(\Om)
\end{equation}
where $\mu_k(\Om)$ is the $k^\text{th}$ non-trivial eigenvalue of 
\begin{equation}
    \label{eq:neumann_ev}
    \begin{dcases}
        - \Delta u = \mu_k(\Om) u &\mbox{ in } \Om, \\
        \partial_n u = 0 &\mbox{ on } \partial \Om.
    \end{dcases}
\end{equation}
Note that by putting $\rho \equiv \frac{1}{\Vol(\Om)}$ in \labelcref{eq:neumann_density}, we can use the eigenvalues of the graph Laplacian to approximate those of \labelcref{eq:neumann_ev} in the sense:
$$
    \mu_k(\Om) = \lim_{n \to \infty} \Vol(\Om) \lambda_k^{(n)}.
$$
From a theoretical point of view, \labelcref{pb:neumann} remains challenging. 
For $k=1$ it is known \cite{weinberger_isoperimetric_1956} that the ball is optimal. 
For $k=2$ it was shown in \cite{bucur_maximization_2018} that the optimal solution is the union of two disjoint, identical balls. However, for all $k \geq 3$ even the existence of optimal domains remains unknown. Considering these incomplete results, some efforts have been put into the numerical approximation of \labelcref{pb:neumann}, e.g. \cite{martinet_numerical_2023}, usually via mesh-based methods. 

In order to embed \labelcref{pb:neumann} into our framework, we first reformulate it as an unconstrained problem; Indeed, it is well known \cite[Eq.~(1.22)]{Henrot} that \labelcref{pb:neumann} is equivalent to
\begin{equation}
    \label{pb:neumann_unconstr}
    \max_{\Om \subset \R^d} \left\{ J(\Om) := \Vol(\Om)^{\frac 2 d}\mu_k(\Om)\right\}.
\end{equation}
Under the assumption that $\mu_k(\Om)$ is a simple eigenvalue and $\Om$ has $C^2$ boundary, the shape derivative (see \labelcref{eq:shape_der}) of $\mu_k(\Om)$ can be expressed as  \cite{henrot_variation_2005}:
\begin{equation}
    \label{eqn:shape_derivative}
    \mu_k'(\Omega, V) = \int_{\partial \Omega} \left(|\nabla u_k|^2 - \mu_k(\Omega) u_k^2\right)(V\cdot n) \dd\Hs
\end{equation}
where $u_k$ is the normalized eigenfunction associated to $\mu_k(\Om)$. 
Using the product rule the shape derivative of \labelcref{pb:neumann_unconstr} equals $ J'(\Omega, V) = \int_{\partial \Omega} f (V\cdot n) \dd \Hs$ where
\begin{equation}    
    \label{eq:sd_neumann}
    f:= \frac{2}{d}\Vol(\Om)^\frac{2-d}{d}\mu_k(\Om) + \Vol(\Om)^\frac{2}{d} \left(|\nabla u_k|^2 - \mu_k(\Om) u_k^2\right)
\end{equation}
Putting $V = f n_\Om$ in the previous expression ensures that advecting $\Om$ by the vector field $V$ will locally increase $J$. 
In order to numerically approximate $V$, the volume $\Vol(\Om)$ and extended normal $n_\Om$ are handled as explained in \cref{subseq:geom}. 
The eigenvalue $\mu_k(\Om)$ is approximated by its graph Laplacian counterpart, as well as the eigenfunction and its gradient for which the values at the vertices are extrapolated to $\R^d$ using a nearest neighbor interpolation.
\begin{remark}[Multiple eigenvalues]
    When $\mu_k(\Om)$ is not simple, it is not necessarily shape differentiable \cite{henrot_variation_2005} which introduces numerical instabilities. Following a similar idea as in \cite{martinet_numerical_2023}, we use a simple regularization method procedure based on the fact that symmetric combinations of eigenvalues are differentiable. When $m$ eigenvalues get close to $\mu_k(\Om)$ during optimization process, we approximate $\mu_k(\Om)$ by a soft minimum depending on a \textit{temperature} parameter $\beta>0$:
    \begin{align*}
        \mu_k(\Om) \approx\frac{\sum_{i=k}^{k+m-1} \mu_i(\Om) e^{-\beta \mu_i(\Om)}}{\sum_{i=k}^{k+m-1} e^{-\beta \mu_i(\Om)}}.
    \end{align*}
\end{remark}

\begin{figure}[!h]
    \centering
    \begin{subcaptionblock}{0.24\textwidth}
        \includegraphics[width=\textwidth]{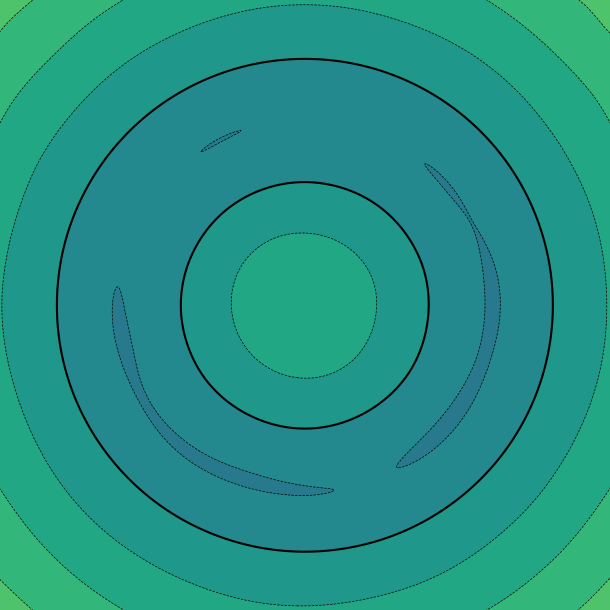}%
        \caption{Initial shape}
    \end{subcaptionblock}
    \begin{subcaptionblock}{0.24\textwidth}
        \includegraphics[width=\textwidth]{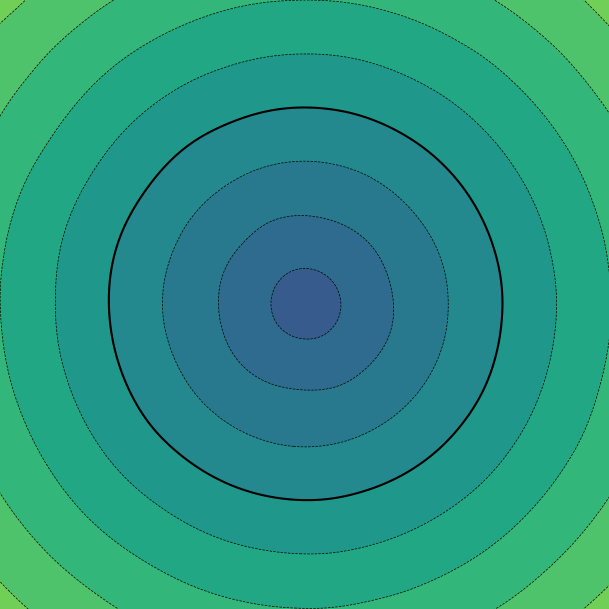}%
        \caption{Optimum for $\mu_1$}
    \end{subcaptionblock}
    \begin{subcaptionblock}{0.24\textwidth}
        \includegraphics[width=\textwidth]{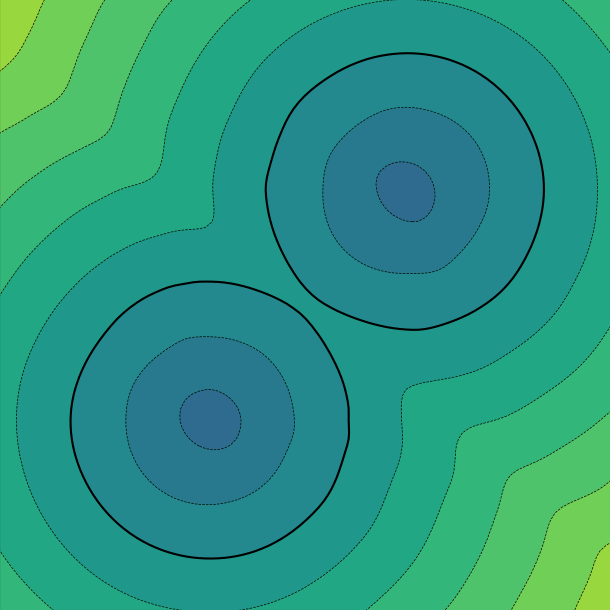}%
        \caption{Optimum for $\mu_2$}
    \end{subcaptionblock}
    \begin{subcaptionblock}{0.24\textwidth}
        \includegraphics[width=\textwidth]{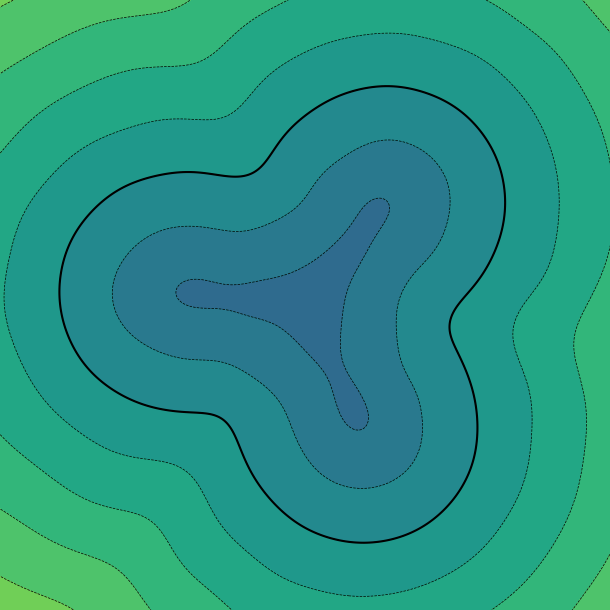}%
        \caption{Optimum for $\mu_3$}
    \end{subcaptionblock}
    \caption{
        \label{fig:neumann_2d}
        2D optimization of Neumann eigenvalues. The bold black lines represent the $0$-level sets.
    }
\end{figure}

\begin{figure}[!h]
    \centering
    \begin{subcaptionblock}{0.24\textwidth}
        \includegraphics[width=\textwidth]{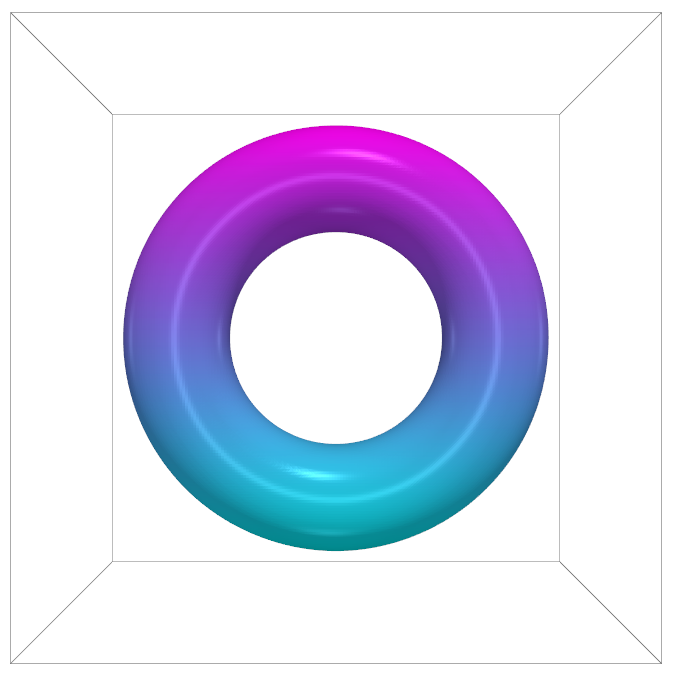}%
        \caption{Initial shape}
    \end{subcaptionblock}
    \begin{subcaptionblock}{0.24\textwidth}
        \includegraphics[width=\textwidth]{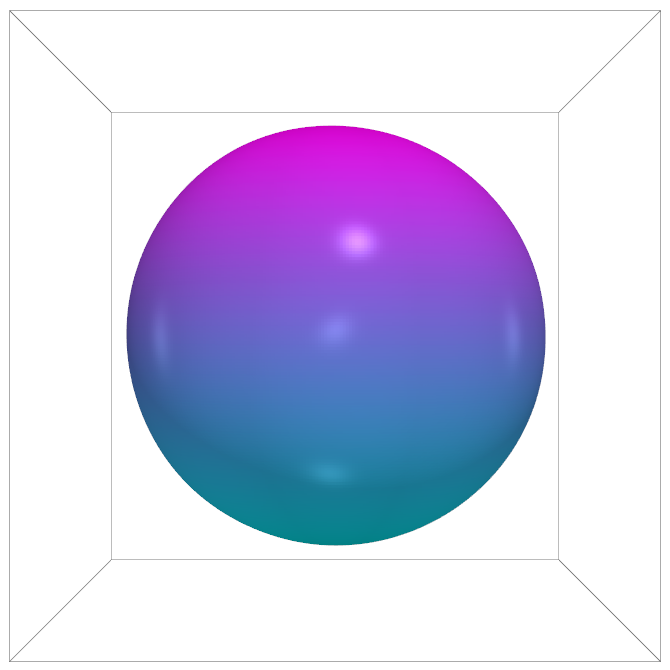}%
        \caption{Optimum for $\mu_1$}
    \end{subcaptionblock}
    \begin{subcaptionblock}{0.24\textwidth}
        \includegraphics[width=\textwidth]{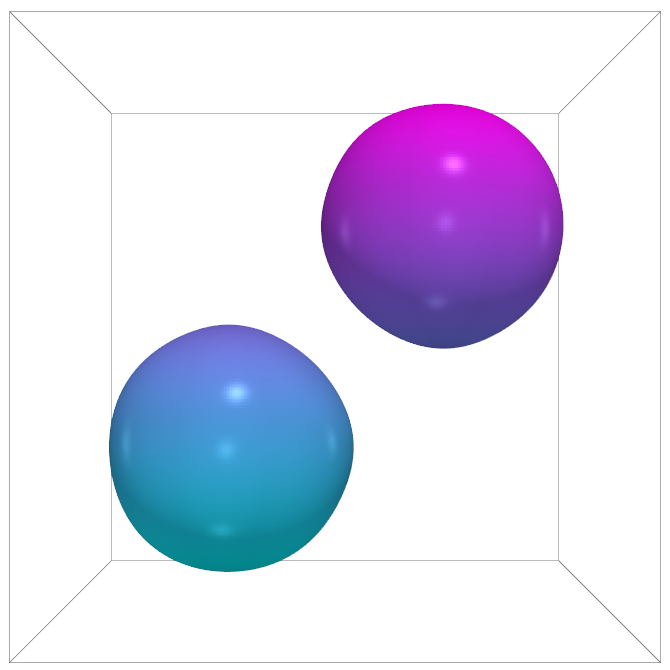}%
        \caption{Optimum for $\mu_2$}
    \end{subcaptionblock}
    \begin{subcaptionblock}{0.24\textwidth}
        \includegraphics[width=\textwidth]{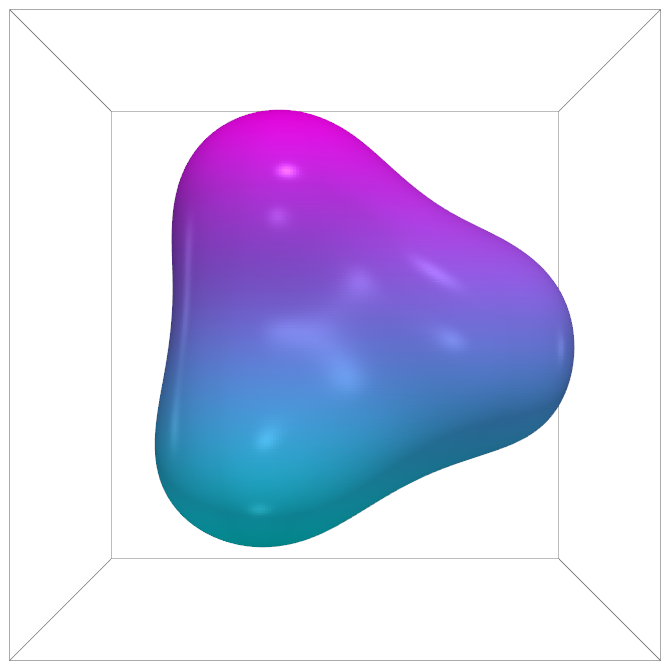}%
        \caption{Optimum for $\mu_3$}
    \end{subcaptionblock}
    \caption{
        \label{fig:neumann_3d}
        3D optimization of Neumann eigenvalues. The depicted surfaces represent the $0$-level set.
    }
\end{figure}

The numerical results of the optimization of \labelcref{pb:neumann} in two and three dimensions are given in \cref{fig:neumann_2d,fig:neumann_3d}. In both cases, we use a 2-layer SIREN network with 80 neurons on each layer, implemented in PyTorch \cite{pytorch}. SIREN networks have already proven to be efficient in representing signed distance functions \cite{sitzmann_implicit_2020}. The results shown here are the ones obtained after $200$ iterations, where each iteration took between $20$ to $40$ seconds to compute on a standard (AMD Ryzen 7 Pro) CPU. We display the results for $k \in \{1,2,3\}$ and, in particular, the algorithm successfully recovers the known optimal shapes for $k=1$ and $k=2$, and the shape for $k=3$ coincides with the one computed in \cite{antunes_numerical_2012}.
We note that our method flawlessly handles the topology changes necessary to get from the initial guess (an annulus) to the optimal shapes and works nicely also in three dimensions.

\subsection{Poisson Problem with Dirichlet Boundary Conditions}

\begin{figure}[!h]
    \centering
    \begin{subcaptionblock}{0.22\textwidth}
        \includegraphics[width=\textwidth]{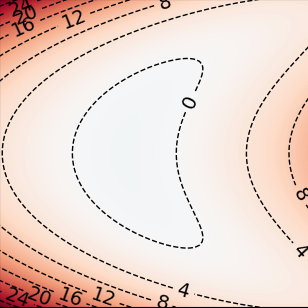}%
        \caption{Function $f$}
    \end{subcaptionblock}
    \begin{subcaptionblock}{0.22\textwidth}
        \includegraphics[width=\textwidth]{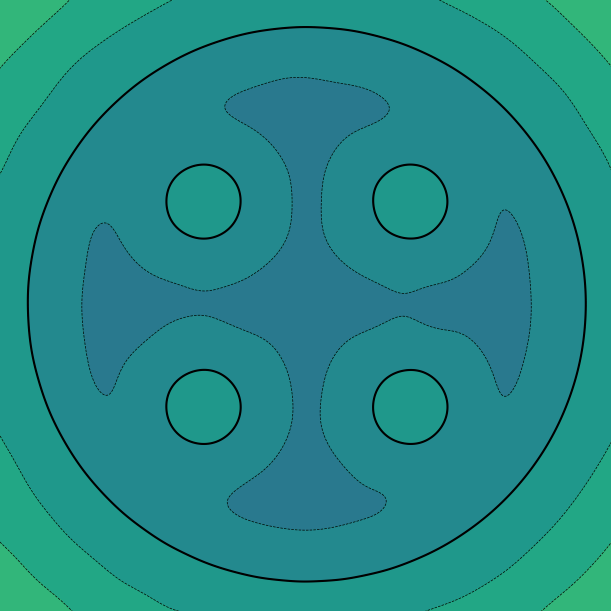}%
        \caption{Initial shape}
    \end{subcaptionblock}
    \begin{subcaptionblock}{0.22\textwidth}
        \includegraphics[width=\textwidth]{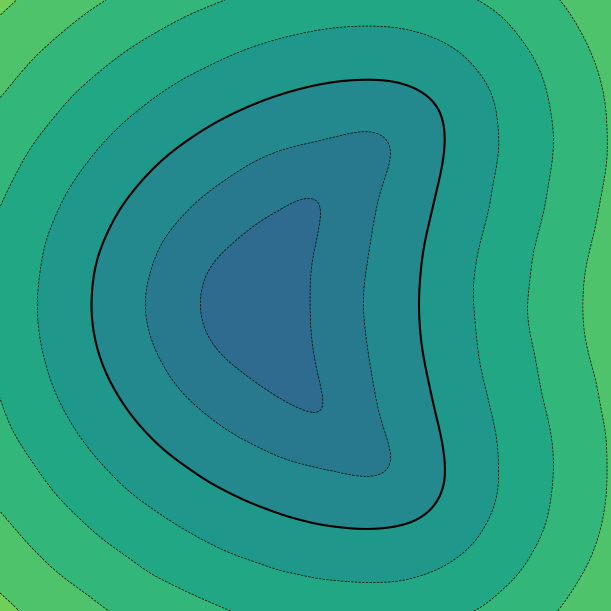}%
        \caption{Optimal shape}
    \end{subcaptionblock}
    \begin{subcaptionblock}{0.29\textwidth}
        \includegraphics[width=\textwidth]{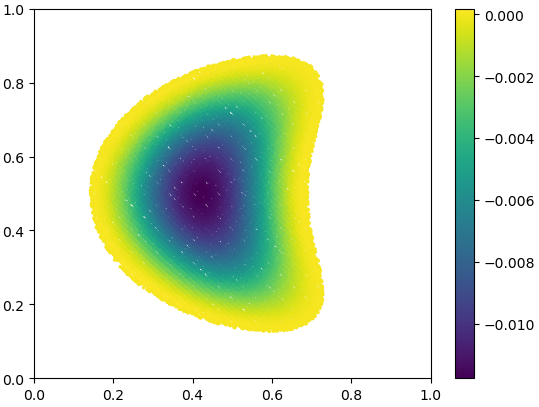}%
        \caption{State function $u$}
    \end{subcaptionblock}
    \caption{
        \label{fig:poisson}
        Results of \eqref{pb:poisson}.
    }
\end{figure}

It turns out that we can also adapt the graph Laplacian to Dirichlet boundary conditions. Let $f \in L^2(\R^d)$ and consider the problem~ \cite{etling2020first}
\begin{equation}
    \label{pb:poisson}
    \min_{\Om \subset \R^d} \int_\Om u\dd x 
\end{equation}
where $u$ is solution to 
\begin{equation} 
    \label{eq:state}
    \begin{dcases}
        - \Delta u = f &\mbox{ in } \Om, \\
        u= 0 &\mbox{ on } \partial \Om. 
    \end{dcases}
\end{equation}
The shape derivative of $J(\Om) := \int_\Om u\dd x$ equals
$
    J'(\Om,V) = \int_{\partial \Om} f(V\cdot n)\dd\Hs
$
where $f := - (\partial_n u)(\partial_n p)$ and
\begin{equation}
    \label{eq:adjoint}
    \begin{dcases}
        \Delta p = 1 &\mbox{ in } \Om, \\
        p = 0 &\mbox{ on } \partial \Om. 
    \end{dcases}
\end{equation}
In order to impose the homogeneous Dirichlet boundary conditions in the graph Laplacian discretization of  \labelcref{eq:state,eq:adjoint}, we use an approach similar to \cite{Jiang2020Jun,calder2023rates}. Set $\Om^\eps := \left\{ x \in \R^d : \d(x,\Om) < \eps \right\}$ and let $x_1,\dots,x_N$ be uniformly i.i.d. samples of $\Om^\eps$ and $L$ be the associated graph Laplacian. We can approximate \labelcref{eq:state} by $\bar L \bar u = \bar f$ where $\bar f \in \R^N$ is such that
$$  
    \bar f_i =
    \begin{dcases}
        f(x_i) & \mbox{ if } x_i \in \Om, \\
        0 &\mbox{ if } x_i \in \Om^\eps\setminus\Om,
    \end{dcases}
$$
and $\bar L$ is the matrix built by replacing the $i^\text{th}$ line of $L$ by the vector $e_i^T$ for all $i$ with $x_i\in\Omega^\eps\setminus\Omega$. 
The function $u$ and its gradient $\nabla u$ can then be approximated from $u^\eps$ as in the previous section and the approximation of $p$ and $\nabla p$ are carried out analogously. The results in two dimensions are presented in \cref{fig:poisson}.

\subsection{Convexity constraints}

\begin{figure}[!h]
    \centering
    \begin{subcaptionblock}{0.32\textwidth}
        \includegraphics[width=\textwidth]{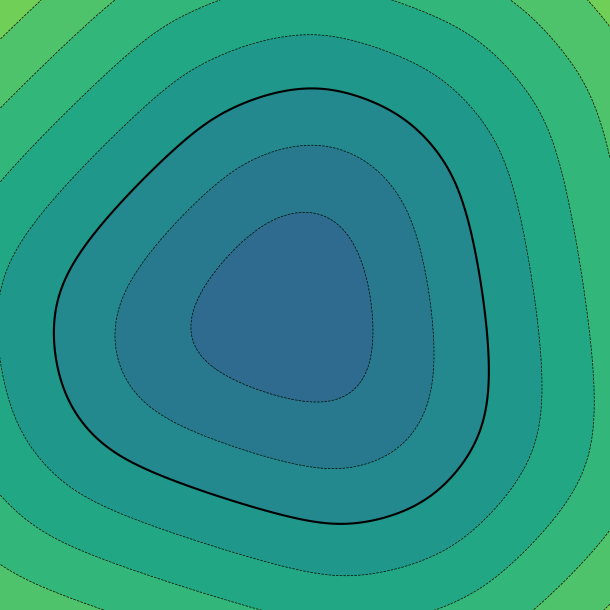}%
        \caption{Optimum for $\mu_2$}
    \end{subcaptionblock}
    \begin{subcaptionblock}{0.32\textwidth}
        \includegraphics[width=\textwidth]{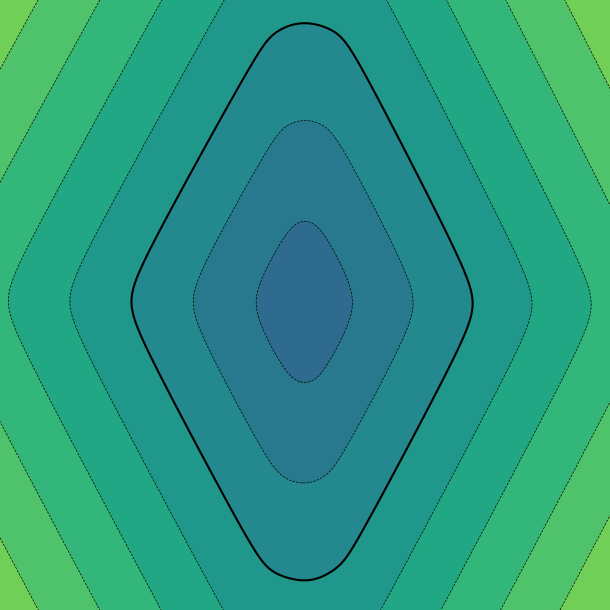}%
        \caption{Optimum for $\mu_3$}
    \end{subcaptionblock}
    \begin{subcaptionblock}{0.32\textwidth}
        \includegraphics[width=\textwidth]{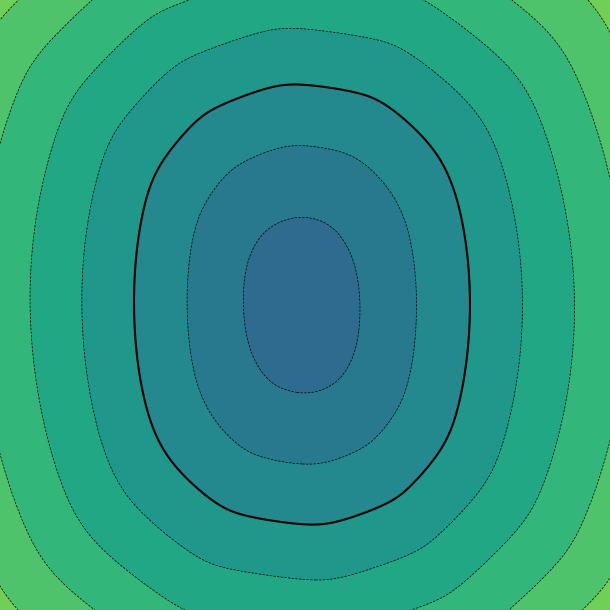}%
        \caption{Optimum for $\mu_4$}
    \end{subcaptionblock}
    \caption{
        \label{fig:convex}
        Numerical solutions to \labelcref{pb:convex}.
    }
\end{figure}

Recently, numerical shape optimization of functionals depending on a convex shape has gained a lot of attention. The usual method to perform numerical shape optimization of convex sets consists in discretizing the support or gauge function in the Fourier or spherical harmonics basis~\cite{antunes_parametric_2022,bogosel_optimization_2024}. Here we show how our framework simply adapts to such shape optimization problems by imposing that the level set function $\phi_\theta$ is convex. For $k\geq 2$, we consider
\begin{equation}
    \label{pb:convex}
    \min_{\Per(\Om)=1} \mu_k(\Om),
\end{equation}
where $\Om \subset \R^d$ is a  convex set.
Let $N\in \N$ , $W_1 \in \R^{N \times d} , W_2 \in \R^{1\times N}, b_1 \in \R^N, b_2 \in \R$. For $\sigma : \R \to \R$ we consider the following 1-layer neural network:
\begin{equation}
    \label{eq:convex_ls}
    \phi_\theta(x) := W_2 \sigma(W_1x+b_1) + b_2.
\end{equation}
If $\sigma$ is convex and all entries of $W_2$ are non-negative, then $\phi_\theta : \R^d \to \R$ is a convex function and the shape $\Om := \{\phi_\theta < 0\}$ is convex. 
As before we can compute the numerical optimizers of \labelcref{pb:convex} for different values of $k$. The results are given in \cref{fig:convex} and should be compared to the ones in \cite{bogosel_optimization_2024}.

\section{Conclusion and Outlook}

In this work we presented a completely mesh-free shape optimization framework that encodes shapes as level sets of a neural network and approximates the occuring PDE using the graph Laplacian. 
We showcased the versatility of our method for three different shape optimization problems. 

Future work will include quantitative comparisons to mesh-based methods, applications to elasticity and surface PDEs, the coupling with other meshless methods like the method of fundamental solutions, and applications to high-dimensional shape optimization problems and data science.

\begin{credits}
\subsubsection{\ackname} 

LB acknowledges funding by the German Ministry of Science and Technology (BMBF) under grant agreement No. 01IS24072A (COMFORT) and by the Deutsche Forschungsgemeinschaft (DFG, German Research Foundation) – project number 544579844 (GeoMAR). 

\subsubsection{\discintname}
The authors have no competing interests to declare that are
relevant to the content of this article.
\end{credits}

%
%
%
\bibliographystyle{splncs04}
\bibliography{bibliography}
%

\end{document}